\documentclass[leqno]{article}

\usepackage{amssymb,eucal,xspace}
\usepackage{amsmath}
\usepackage{amsthm}

\swapnumbers
\newtheorem{num}{}[section] 

\newcommand{\beq}{\addtocounter{num}{1}\begin{equation}}
\newcommand{\alku}{\begin{num}\it }
\newcommand{\loppu}{\end{num}}

\newcommand{\spa}{{\rm span}\,}
\newcommand{\cspa}{\overline{\rm span}\,}
\newcommand{\id}{{\rm id}}
\def\ppav]#1]{\left]#1\right]}
\newcommand{\pf}{\it Proof. \rm }
\newcommand{\imp}{$\Rightarrow$\xspace}
\newcommand{\cd}{{\rm cd}\, A}
\newcommand{\sub}{\subset}
\newcommand{\tyh}{\emptyset}
\newcommand{\sm}{\setminus}
\newcommand{\lau}{{\bf Theorem.} }
\newcommand{\lem}{{\bf Lemma.} }
\newcommand{\cor}{{\bf Corollary.} }
\newcommand{\iy}{\infty}
\newcommand{\rn}{{\sf R}^n}
\newcommand{\eni}{$\varepsilon$-near\-isometry\xspace}

\begin{document}






\vspace*{5cm}
\begin{center}
\Large
{\bf Isometric approximation property of unbounded sets}

\bigskip

\normalsize
Jussi V\"ais\"al\"a
\end{center}

\bigskip

\begin{abstract}
We give a necessary and sufficient quantitative geometric condition for an unbounded
set $A \sub\rn$ to have the following property with a given $c >0$: For every
$\varepsilon \ge 0$ and for every map $f\colon A\to\rn$ such that $\big| |fx-fy| -
|x-y|\big| \le \varepsilon$ for all $x,y\in A$, there is an isometry $T\colon  A
\to\rn$ such that $|Tx - fx| \le c\varepsilon$ for all $x \in A$.
\medskip

\end{abstract}

\section{Introduction}

\alku
\label{1.1}
Near\-iso\-metries. \rm
Let $E$ and $F$ be real Hilbert spaces and let $ A \sub E$. A map
$f\colon A \to F$ is said to be a \textit{near\-isometry} if there is
$\varepsilon \ge 0$ such that 
\[
|x-y| - \varepsilon \le |fx-fy| \le | x-y|+ \varepsilon
\]
for all $x,y \in A$. More precisely, such a map is an
$\varepsilon$-\textit{near\-isometry}. In the literature, the
$\varepsilon$-near\-isometries are often called $\varepsilon$-isometries.

We study the stability question: How well can an \eni be approximated by a true
iso\-metry? The following fundamental result was proved by D.H. Hyers and S.M. Ulam
\cite{HU} in 1945: If $f\colon E \to F$ is a surjective \eni, then there is a surjective
iso\-metry $T\colon E \to F$ with
\[
d(T,f) = \sup \{  |Tx-fx|: x \in E\} \le 10\varepsilon.
\]
The bound $10\varepsilon$ has  later been improved to $\sqrt 2 \varepsilon$; see
\cite[1.2]{Se}. For a Banach space version, see \cite[15.2]{BL}. In the case $E=F=\rn$
the result holds without the surjectivity condition, as proved by R. Bhatia
and P. \v Semrl \cite[Th.1]{BS}. A survey of near\-isometries is given in \cite{survey}.

The case where $ A \sub E=F=\rn$ 
 has been studied in \cite{ATV} and in
\cite{solar}. For $c > 0$ we say that a set $A \sub \rn$ has the
$c$-\textit{isometric approximation property}, abbreviated $c$-IAP, if for each
$\varepsilon \ge 0$ and for each \eni
$f\colon A \to \rn$ there is an iso\-metry $T \colon \rn\to\rn$ such that
\[
d(T,f) = \sup \{  |Tx-fx|: x \in A\} \le c\varepsilon.
\]
If $A$ has the $c$-IAP and contains at least two points, then clearly $c \ge 1/2$.

The whole space $\rn$ has the $\sqrt 2$-IAP. In \cite[2.5]{solar} we gave a quantitative
geometric characterization for \textit{bounded} subsets of $\rn$ with the $c$-IAP  
 in terms of the so-called $c$-solar systems. For sets without
isolated points, it can be expressed in the following simple form. For a unit vector
$u\in\rn$ let $\pi_u\colon \rn\to {\sf R}$ be the projection $\pi_ux = x \cdot u$. The
\textit{thickness} of a bounded set $A \sub \rn$ is the number
\[
\theta(A) = \inf \{  d(\pi_u A): |u| = 1\},
\]
where $d$ denotes diameter. Then $A$ has the $c$-IAP if and only if
$\theta(A) \ge d(A)/c'$, where $c$ and $c'$ depend only on each other and on $n$.
In particular, a ball has the $c$-IAP with $c=c(n)$; this was essentially
proved by F. John \cite{Jo} already in 1961. From a recent result of E. Matou\v
skov\'a \cite[4.1]{Ma} it follows that $c(n)$ cannot be chosen to be independent of $n$.

The motivation for the present paper was to find a characterization for the
\textit{unbounded} subsets of $\rn$ with the $c$-IAP. It turns out that the answer is
simpler than in the case of bounded sets, and it is given in \ref{main}. 
We associate to each unbounded set
$A \sub\rn$ a  number $\mu(A) \in [0,1]$ and show that if $A$ has the $c$-IAP, then
$\mu(A)
\ge 1/c' > 0$ with $c' = 17c$. Conversely, if $\mu(A) \ge 1/c'$, then $A$ has the
$c$-IAP with $c =  \sqrt 2 c'$. The constant $\sqrt 2$ is the best possible. 

Contrary to
the bounded case, these estimates do not depend on the dimension $n$. In fact, several
arguments are valid in arbitrary Hilbert spaces,  but the most satisfactory formulation
is obtained  in the euclidean $n$-space $\rn$.
\loppu

\alku
\label{summary}
Summary of the paper. \rm
We start by introducing in Section 2   the number $\mu(A)$ and other auxiliary
concepts needed in the proof. We can then state the main theorem  \ref{main} and various
related remarks. In Section 3 we show that $c$-IAP implies the property
$\mu(A)
\ge 1/c'$. The converse is proved in Section 4, where we give a more
general result valid for all Hilbert spaces,  which may have independent interest. In
Section 5 we give a short elementary proof and an improvement of the original
Hyers-Ulam theorem. A remark on the Banach space case is given in Section 6.
\loppu

\section{Formulation of the main result} 

\alku
\label{2.1}
Notation. \rm
Let $E$ be a real Hilbert space.
We let $x \cdot y$ denote the  inner product of vectors $x,y \in E$, and the
norm of $x$ is $|x| = \sqrt{x \cdot x}$. Open and closed balls with center $x$ and
radius $r$ are written as $B(x,r)$ and $\bar B(x,r)$, respectively, and we abbreviate
$B(r) = B(0,r),\ \bar B(r) = \bar B(0,r)$. The unit sphere is $S(1) = \partial
B(1)$. The diameter of a set $A\sub E$ is $d(A)$, and the distance between nonempty
sets
$A,B
\sub E$ is $d(A,B)$. The standard basis of
$\rn$ is
$e_1,\dots,e_n$. The central projection
$p\colon
 E \sm\{  0\} \to S(1)$ is defined by
\[
px = x/|x|.
\]
An expression like $ab/cde$ means $(ab)/(cde)$.
To simplify expressions we often omit parentheses writing  $fx = f(x)$ etc.
\loppu

\alku
\label{2.2}
Definitions. \rm
Let $A \sub E$ be an unbounded set and let $u \in E$ be a unit vector. We say that a
sequence $(x_j)$ in $A$ \textit{converges directionally} to $u$ if $|x_j| \to\iy$ and
$px_j \to u$. A unit vector $u$ is a \textit{cluster direction} of $A$ if there is a
sequence in $A$ converging directionally to $u$.
 We let $\cd$ denote the set of all cluster directions of $A$. Equivalently,
\[
\cd = \bigcap \{  {\rm cl}\, p[A \sm B(R)]: R > 0\}.
\]
 Obviously ${\rm cd}\, (A+b) = \cd$ for
each $b \in E$. If $\dim E < \iy$, the set $\cd$ is  compact and nonempty.

Let $ X \sub S(1)$. For $e \in S(1)$ we set
\[
\sigma(e,X) = \sup \{  |u \cdot e|: u \in X\},\quad \mu_1(X) = \inf \{  \sigma(e,X): |e|
= 1\}
\]
with the convention $\mu_1(\tyh) = 0$. Thus $\mu_1(X)$ is small if and only if $X$ lies
in a narrow neighborhood of a hyperplane through the origin. In fact, $\mu_1(X) =
\theta(X
\cup (-X))/2$. 

For an unbounded set $A \sub E$ we set
\[
\mu(A) = \mu_1(\cd).
\]
For alternative definitions of $\mu(A)$ for $A \sub\rn$, see \ref{alt}.

We can now state the main result of the paper:
\loppu

\alku
\label{main}
{\bf Main theorem.}
For an unbounded set $A \sub \rn$, the following conditions are quantitatively
equivalent:

$(1)$ $A$ has the $c$-{\rm IAP}.

$(2)$  $\mu(A) \ge 1/c'$.

\noindent More precisely, $(1)$ implies $(2)$ with $c' = 17c$, and $(2)$ implies $(1)$
with $c =  \sqrt 2 c'$. The constant $\sqrt 2$ is the best possible.
\loppu

Part (1) \imp (2) will be proved in Section 3, and the converse part in Section 4.
Observe that the bounds in \ref{main} do not depend on $n$.

\alku
\label{2.4} Remarks. \rm
1. Always $0 \le \mu(A) \le 1$.

2. The number $\mu(A)$ is invariant under translations: $\mu(A + b) = \mu(A)$ for each
$b \in\rn$.

3. If $A \sub A'$, then $\mu(A) \le \mu(A')$.

4. If $A$ contains a half space, then $\mu(A) = 1$.

5. For $e \in S(1)$ and $0 \le \alpha \le \pi/2$, let $C(e,\alpha)$ be the cone 
$
 \{  x \in E: x \cdot e \ge |x| \cos \alpha\}.
$
Then $\mu(C(e,\alpha)) = \sin \alpha$.

6. Since $\mu(\rn) = 1$, the Bhatia-\v Semrl result mentioned in \ref{1.1} follows
from
\ref{main} as a corollary.

\loppu

\alku
\label{2.4+1}
Sharpness. \rm
The bound $c = \sqrt 2c'$ in \ref{main} is the best possible. In fact, for each $n$	
there is a 1-near\-isometry $f\colon \rn\to\rn$ such that $d(T,f)$ is at least the
Jung constant $J(\rn) = \sqrt{2n/(n+1)}$ of $\rn$ for every iso\-metry $T\colon
\rn\to\rn$; see
\cite[1.11]{HV}. Hence the whole space $\rn$ does not have the $c$-IAP for any $c <
J(\rn)$. Since $\mu(\rn) = 1$, the number $\sqrt 2$ cannot be replaced by any smaller
universal constant.

On the other hand, the bound $17c$ is presumably far from optimal. We get a lower
estimate for this bound by considering the set $A = {\sf R}$, which has the $1$-IAP
in ${\sf R}$. Since $\mu(A) = 1$, Theorem
\ref{main} is not true if the number 17 is replaced by any universal constant less
than 1.
\loppu

\alku
\label{2.4+2}
Uniqueness. \rm 
Suppose that $A \sub\rn$ is an unbounded set with the $c$-IAP.
If $f\colon A \to\rn$ is an \eni, then the iso\-metry $T$ with $d(T,f) \le
c\varepsilon$ given by the $c$-IAP is uniquely determined up to translation. To prove
this, it suffices to show that there is at most one linear iso\-metry $U\colon
\rn\to\rn$ with $d(U,f) \le K < \iy$. Assume that $U$ and $U'$ are such maps. Let $e
\in \cd$. Since $\mu(A) > 0$, the set $\cd$ spans $\rn$. Hence it suffices to show
that $Ue = U'e$.

Choose a sequence $(x_j)$ in $A$ converging directionally to $e$. Since
$|Ux-U'x| \le 2K$ for all $x \in A$, we get $|U(px_j) - U'(px_j)| \le 2K/|x_j|$.
Letting $j \to\iy$ yields $Ue - U'e = 0$ as desired.

\loppu

\alku
\label{alt}
Alternative definitions for $\mu(A)$. \rm
 The following considerations may clarify the
meaning of $\mu(A)$, but they are not actually needed in the paper. 

Let $E$ be a Hilbert space, let $e \in E$ be a unit vector, and let $0 \le \alpha \le
\pi/2$. We let $D(e,\alpha)$ denote
 the open double cone with
axis ${\rm span}\, e$, vertex 0 and central angle $\alpha$:
\[
D(e,\alpha) = \{  x \in\rn: |x \cdot e| > |x| \cos \alpha\} = \bigcup \{  B(te,|t|
\sin \alpha): t \in {\sf R}\}.
\]
For an unbounded set $A \sub E$ we write
\begin{align*}
\varphi(A) &= \sup \{  \alpha: A \cap D(e,\alpha)  \text{ is bounded for some } e \in
S(1)\},\\
 \tau(A) &= \sup_{|e|=1} \liminf_{|t|\to\iy} d(te,A)/|t|.
\end{align*}

\alku
\label{2.8}
{\bf Proposition.} Always $\tau(A) = \sin \varphi(A),\ \mu(A) \le \cos \varphi(A)$. If
$\dim E < \iy$, then 
$
\mu(A) = \cos \varphi(A) = \sqrt{1-\tau(A)^2}.
$
\loppu

We omit the easy proof. The inequality $\mu(A) \ge \cos \varphi(A)$ for $E=\rn$ follows
from
\ref{3.1} below.
\loppu

\section{Proof for (1) \imp (2)}

To obtain the part (1) \imp (2) of \ref{main} we define an auxiliary map and prove
several lemmas.

\alku
\label{3.1}
\lem
Suppose that $L \sub\rn$ is a linear subspace, and let $P\colon \rn\to L$ and $P'\colon
\rn\to L^\perp$ be the orthogonal projections.
Suppose also that $A \sub \rn$ is an unbounded set such that that 
$|P'u| < q
\le 1$ for all $u \in \cd$. Then the set $A \cap \{  x : |P'x| \ge q|x|\}$ is
bounded. 
\loppu

\pf
If the lemma is false, there is a sequence $(x_j)$ in $A$ converging directionally to
$u \in\cd$ such that $|P'x_j| \ge q |x_j|$ for all $j$. Dividing by $|x_j|$ and
letting $j\to\iy$ yields $|P'u| \ge q$, a contradiction. 
  $\square$

\alku
\label{3.1+1}
Remarks. \rm
1. The bounded set in \ref{3.1} can also be written as $A \cap \{  x: |P'x| \ge r|Px|\}$,
where $r=q/\sqrt{1-q^2}$.

2. Lemma \ref{3.1} does not hold in infinite-dimensional Hilbert spaces, because  an
unbounded set need not have any cluster directions. This is the main reason why we
consider in this section only the case $E=\rn$.

\loppu

\alku
\label{3.2}
An auxiliary map. \rm
In \ref{3.2}--\ref{3.9} we assume that $n \ge 2$ and write $\rn = {\sf R}^{n-1} \times
{\sf R}$. Given a number $M \ge 1$
 we define a map $g\colon \rn\to\rn$ by
\begin{equation*}
g(x,t) = 
  \begin{cases}
(x, t + \sqrt{|x|}) &\text{ for } |x| \le M,\\
(x, t + \sqrt M)  &\text{ for } |x| \ge M.
  \end{cases}
\end{equation*}
We want to show that $g$ is a 1-near\-isometry in a certain set. For that purpose,
let $z = (x,t),\ z' = (x',t') \in\rn$. Set
\[
d = |z'-z|,\ D = |gz' - gz|, \ h = |x'-x|.
\]
\loppu

\alku
\label{3.3}
\lem
If $|x| \le |x'| \le M,\ r=1/8\sqrt M,\ |t| \le r|x'|$ and $|t'| \le 3r|x'|$, then
$|D-d|
\le 1$.
\loppu

\proof
A direct computation gives $D^2 = d^2 + a$, where 
\[
a = |x| + |x'| - 2\sqrt{|x||x'|} + 2(t'-t) (\sqrt{|x'|} - \sqrt{|x|}).
\]
Since $|x| \le |x'|$ and
\[
\big|\sqrt{|x'|} - \sqrt{|x|}\big| = \frac{\big| |x'| -|x|\big|}{\sqrt{|x'|} +
\sqrt{|x|}}
\le
\frac{h}{\sqrt{|x'|}},
\]
we get
\[
a \le |x'| - |x| + 2(|t'|+|t|)h/\sqrt{|x'|}
\le h + 8rh\sqrt M = 2h.
\]
Furthermore,
\[
a \ge -2(|t'|+|t|)(\sqrt{|x'|} - \sqrt{|x|}) \ge -h,
\]
and hence $|a| \le 2h$. Since $D \ge h$ and $d \ge h$, this implies that
\[
|D-d| = \frac{|D^2 - d^2|}{D+d} \le \frac{|a|}{2h} \le 1.\quad  \square
\]

\alku
\label{3.4}
Notation. \rm
For $r \ge 0$ we set
\[
A_r = \{  (x,t) \in\rn: |t| \le r|x|\}.
\]
Thus $A_r = \rn \sm D(e_n,\alpha)$, where $D(e_n,\alpha)$ is the double cone defined in
\ref{alt} and $\cot \alpha = r$. Setting $q = r/\sqrt{1 + r^2}$ we can also write
\beq
\label{Ar}
A_r = \{  z \in \rn: |z \cdot e_n| \le q|z|		\}.
\end{equation}
Observe that $r = q/\sqrt{1-q^2}$.
\loppu

\alku
\label{3.4+1}
\lem
Suppose that $M > 0$ and that $x,x' \in\rn$ with $|x| < M,\ |x'| > M$. Let $x''$ be the
point of the line segment $[x,x']$ with $|x''| = M$. Then 
\[
|x'||x-x''| \le 2M|x-x'|.
\]
\loppu
\pf
Set $s = |x'||x-x''|/|x-x'|$. We consider three cases.

\textit{Case} 1. $|x'| \le 2M$. Since $|x-x''| \le |x-x'|$, we have $s \le |x'| \le 2M$.

\textit{Case} 2. $|x'| \le |x-x' |$. Now $s \le |x-x''| \le |x|+|x''| \le 2M$.

\textit{Case} 3. $|x'| \ge 2M$ and $|x-x'| \le |x'|$. Since $|x''-x'| \ge |x'| - M$, we
have $|x-x''| = |x-x'| - |x''-x'| \le M$.
  Thus
\[
s \le \frac{M|x'|}{|x'| - |x|} = \frac{M}{1 - |x|/|x'|} \le 2M. \quad \square
\]

\alku
\label{3.5}
\lem
Suppose that $M \ge 1$ and that $r = 1/8\sqrt M$. Then the map $g|A_r$ is a
$1$-near\-isometry.
\loppu

\pf
Let $z = (x,t),\ z' = (x',t') \in A_r$ with $|x| \le |x'|$, and let $d$ and $D$ be as
in \ref{3.2}. We must show that $|D-d| \le 1$. If $|x'| \le M$, this follows from
\ref{3.3}. If $|x| \ge M$, then $|D-d| = 0$. It remains to consider the case $|x| \le
M \le |x'|$.

There is a point $z'' = (x'',t'')$ of the line segment $[z,z']$ with $|x''| = M$.
Since
$|t| \le r|x| \le rM$, we obtain by \ref{3.4+1}
\[
|t''| = \Big| \frac{|x'-x''|}{|x'-x|} t + \frac{|x-x''|}{|x'-x|} t' \Big| \le Mr + 2Mr =
3r|x''|.
\]
 Hence $|gz''-gz|\le |z''-z|+1$ by \ref{3.3}.
Consequently,
\[
D \le |gz'-gz''| + |gz''-gz| \le |z'-z''| + |z''-z| +1 =d+1.
\]
It remains to show that $D \ge d-1$. Computing as in \ref{3.3} we get $D^2 = d^2 + a$
where
\[
a \ge 2(t'-t)(\sqrt M - \sqrt{|x|}) \ge -2r(|x'|+|x|)(\sqrt M - \sqrt{|x|}).
\]
Arguing as in \ref{3.3} we see that it suffices to show that $a \ge -2h = -2|x'-x|$. We
show that
$a \ge -h$.

\textit{Case} 1. $M \le h$. Now $|x'|+|x| \le |x'-x|+2|x| \le h+2M \le 3h$, and hence
\[
a \ge -2r \cdot 3h\sqrt M = -6rh/8r \ge -h.
\]

\textit{Case} 2. $M \ge h$. Now $|x'|+|x| \le h+2M \le 3M$. Since $M-|x| \le h$, we get
\[
a \ge -6rM \frac{M-|x|}{\sqrt M + \sqrt{|x|}} \ge -6rh\sqrt M = -6h/8 \ge -h.\quad 
\square
\]

\alku
\label{3.6}
Notation \rm
We fix a number $\lambda,\ 0 < \lambda \le 1$; later we let $\lambda\to 0$. Let $r$ be
as in
\ref{3.5} and set $A'(r,\lambda) = A_r \cup \bar B(\lambda)$. Define $f\colon 
A'(r,\lambda) \to\rn$ by
\begin{equation*}
fx =
\begin{cases}
gx & \text{ for } x \in A_r,\\
0 & \text{ for } x \in \bar B(\lambda) \sm A_r.
\end{cases}
\end{equation*}
\loppu

\alku
\label{3.7}
\lem
The map $f\colon A'(r,\lambda)\to\rn$ is a $(1+\lambda)$-near\-isometry.
\loppu

\proof
Let $z,z' \in A'(r,\lambda)$ and set $d = |z'-z|,\ D = |fz'-fz|$. If $z,z' \in
A_r$, then $|D-d| \le 1$ by \ref{3.5}. If $z,z' \in  \bar B(\lambda)\sm A_r$, then
$|D-d| =d \le 2\lambda \le 1+\lambda$. Finally, let $z \in \bar B(\lambda) \sm A_r,\ z'
\in A_r$. Since $f|A_r$ is a 1-near\-isometry by \ref{3.5} and since $f(0)=0$, we
have $\big| |fz'| - |z'|\big| \le 1$, and hence
\[
|D-d| = \big| |fz'| - |z-z'|\big| \le 1+|z| \le 1+\lambda.\quad\square
\]

\alku
\label{3.9}
\lem
Let $f\colon A'(r,\lambda)\to\rn$ be as in {\rm \ref{3.6}}, let $A \sub A'(r,\lambda)$
be unbounded, and let $T \colon \rn\to\rn$ be a linear isometry with $d(T,f|A) < \iy$.
Then $T|\cd = \id$.
\loppu

\proof
Let $u \in \cd$. Set $K = d(T,f|A)$ and choose a sequence $(z_j)$ in $A$ such that
$|z_j|\to\iy$ and $ u_j = pz_j \to u$. For large $j$ we have $fz_j = z_j + w$ with $w =
\sqrt M e_n$. Then $|Tz_j -z_j -w| \le K$, and hence
\[
\big| Tu_j - u_j - w/|z_j|\big| \le K/|z_j|.
\]
As $j \to\iy$, this implies $Tu=u$.  $\square$

\medskip
We turn to the proof of part (1) \imp (2) of Theorem \ref{main}.

\alku
\label{3.10}
\lau
If $A \sub\rn$ is an unbounded set with the $c$-{\rm IAP}, then $\mu(A) \ge 1/17c$.
\loppu

\proof
The theorem is trivially true for $n=1$, since then $\mu(A) = 1$ for each unbounded set
$A \sub {\sf R}$. Assume that $n \ge$ 2. Set $X = \cd$ and $L = \spa X$. We first show
that
$L=\rn$. Assume that
$L\ne\rn$. We may assume that $L\sub {\sf R}^{n-1}$. Let $P\colon \rn\to L$ and
$P'\colon \rn\to L^\perp$ be the orthogonal projections. Set $M = 25c^2$ and $r =
1/8\sqrt M = 1/40c$. The set $Q = \{  x\in A: |P'x|> r|Px|\}$ is bounded by \ref{3.1}.
Choose a number $R> 0$ with $Q \sub \bar B(R)$. Since the property $c$-IAP is a
similarity invariant by
\cite[3.1]{solar}, we may replace $A$ by $A/R$ and thus assume that $Q \sub \bar B(1)$.

Choose $z \in A$ with $|Pz| \ge M$. Write $P'z = te$ with $t \ge 0,\ |e|=1$. Using an
auxiliary orthogonal map keeping $L$ fixed we may assume that $e=e_n$. Then $z = x +
te_n$ with $x = Pz \in {\sf R}^{n-1}$.

Let $f\colon  A'(r,1) \to \rn$ be as in \ref{3.6}. By \ref{3.7} this map is a
2-near\-isometry. Since $A \sub A'(r,1)$ and since $A$ has the $c$-IAP, there
is an iso\-metry $T \colon \rn\to\rn$ with $d(T,f|A) \le 2c$. Since $f(0)=0$, there is
a linear isometry $U\colon \rn\to\rn$ with $d(U,f|A) \le 4c$. By \ref{3.9} we have
$U|L=\id$. Since
\[
Uz = U(x+te_n) = x+tUe_n,\quad fz = x + (t+\sqrt M)e_n,
\] 
we obtain 
\[
4c \ge |fz-Uz| = |(t+\sqrt M) e_n - tUe_n| \ge |(t + \sqrt M)e_n| - |tUe_n|
 = \sqrt
M = 5c,
\]
a contradiction. We have proved that $\spa X = \rn$.

Set $q = 1/17c$ and $r = q/\sqrt{1-q^2}$. Since $c \ge 1/2$ (see \ref{1.1}), we have
$r \le 1/8$ and $M = 1/64r^2 \ge 1$. Assume that $\mu(A) < q$.  Choose $e \in S(1)$
such that $|u \cdot e| \le q$ for all $u \in X$. We may assume that $e = e_n$.
By \ref{3.1} and \ref{3.1+1}.1, the set $A \sm A_r$ is bounded; say $A \sm A_r \sub \bar
B(R)$. Arguing as above, we may replace $A$ by $\lambda A/R$; then $A \sm A_r \sub
\bar B(\lambda)$. Define the $(1+\lambda)$-nearisometry $f\colon  A'(r,\lambda)
\to\rn$ as in \ref{3.6}. As above, there is a linear isometry $U\colon \rn\to\rn$
with $d(U,f|A) \le 2c(1+\lambda) $. Since $\spa X = \rn$, we have $U = \id$ by
\ref{3.9}. 

Choose a point $z = (x,t) \in A$
 with $|x| \ge M$. By the definition \ref{3.2} of $g$, we have
$
|fz-z| = |gz-z| = \sqrt M = 1/8r.
$
On the other hand, $|fz-z| = |fz-Uz| \le 2c(1+\lambda)$. Hence $2c(1+\lambda) \ge
1/8r$. As $\lambda\to 0$, this yields $r \ge 1/16c$. Since $r \le 1/8$, this implies
that $q = r/\sqrt{1+r^2} \ge 1/2c\sqrt{65} > 1/17c$, a contradiction. 
$\square$

\section{Proof for (2) \imp (1)}

\alku
\label{4.1}
Outline of the proof. \rm
We consider the more general problem where $E$ and $F$ are arbitrary Hilbert spaces
and
$A$ is an unbounded subset of $E$ with $\mu(A) \ge 1/c'$. Suppose that $f\colon A \to F$
is an \eni. We shall show that there is an iso\-metry $T\colon E \to F$ such that
$d(T,Pf) \le \sqrt 2\varepsilon$ where $P$ is the orthogonal projection of $F$ onto
$TE$. The part (2) \imp (1) of \ref{main} is then an immediate corollary.

We normalize the situation by the conditions $0 \in A$ and $f(0)=0$.
Set $X = \cd$. We define a map $\varphi\colon X \to S(1) = S_F(1)$ as follows. Let $u \in
X$. Choose a sequence $(x_j)$ in $A$ converging directionally to $u$. We show
that the sequence $(pfx_j)$ converges to a point $u' \in S(1)$ and that $u'$ is
independent of the choice of the sequence $(x_j)$. Setting $\varphi u = u'$ we obtain
the map $\varphi$. This map is an isometry, and it extends to a linear isometry
$U\colon E\to F$. We show that $d(U,Pf) \le 2c'\varepsilon$, where $P\colon F \to UE$
is the projection. Combining
$U$ and a translation we obtain the desired isometry $T$ with $TE = UE$ and  $d(T,Pf)
\le
\sqrt 2 c'\varepsilon$.
\loppu

\alku
\label{4.not} 
Notation. \rm
 In this section $E$ and $F$ will always denote real Hilbert spaces of
dimension at least one.

\loppu

\alku
\label{4.2}
\lem
Suppose that $\varepsilon > 0$ and that $x,y \in E \sm B(2\varepsilon)$. Suppose also
that $f\colon \{  0,x,y\}
\to F$ is an \eni with $f(0) = 0$. Then
\[
|pfx - pfy|^2 \le 4 |px-py|^2 +  24\varepsilon/|x| + 24 \varepsilon/|y| +
12\varepsilon^2/|x||y|.
\]
\loppu

\proof
Set
\[
s = |x|, \ t = |y|,\ r = |x-y|,\ s' = |fx|,\ t' = |fy|,\ r' = |fx-fy|.
\]
Since
\[
|a||b||pa-pb|^2 = |a-b|^2 - (|a|-|b|)^2
\]
for all $a,b \in E \sm\{  0\}$, and since
$|t-s| \le r,\ |t'-s'| \le r'$, we obtain
\begin{align*}
& s't'|pfx-pfy|^2 - st|px-py|^2\\
&= (r+r')(r'-r) + (t+t'-s-s')(t-t'+s'-s)\\
&\le (r+r')|r'-r| +(|t-s| + |t'-s'|)(|t'-t| + |s'-s|)\\
& \le (r+r')(|r'-r| + |t'-t| + |s'-s|) \le (2r+\varepsilon)3\varepsilon.
\end{align*}
Since $r \le s+t,\ 
s' \ge s-\varepsilon \ge s/2$ and $ t' \ge t-\varepsilon \ge t/2
$,
this implies the lemma.\   $\square$

\alku
\label{4.5}
\lem
Suppose that $A \sub E$ is unbounded, that $f\colon A \to F$ is a
near\-isometry with $f(0)=0$, and that $u \in\cd$. Choose a sequence $(x_j)$ in $A$
converging directionally to $u$. Then the sequence
$(pfx_j)$ converges to a point $u' \in S(1)$, and $u'$ is independent of the choice of
the sequence $(x_j)$.
\loppu

\proof 
Since the sequence $(px_j)$ is Cauchy and since $|x_j| \to\iy$, it follows from
\ref{4.2} that the sequence $(pfx_j)$ is Cauchy and hence convergent. If $(y_j)$ is
another sequence in $A$ converging directionally to $u$, then so is the sequence $(z_j)
= (x_1,y_1,x_2,y_2,\dots)$. Since the sequence $(pfz_j)$ is convergent, the
subsequences $(pfx_j)$ and $(pfy_j)$ converge to the same limit.  $\square$

\medskip
We next prove an elementary but useful inequality.

\alku
\label{M}
\lem
Suppose that 
$x,y \in E$ and that $f\colon  \{  0,x,y\} \to F$ is an \eni with $f(0)=0$. Then
\[
|fx \cdot fy - x \cdot y| \le 2\varepsilon(|x| + |y| + \varepsilon).
\]
\loppu

\proof
Since
\[
2x \cdot y = |x|^2 + |y|^2 - |x-y|^2,\quad 2fx \cdot fy = |fx|^2 + |fy|^2 - |fx-fy|^2,
\]
we get
\begin{align*}
2|fx \cdot fy - x \cdot y| &\le \big||fx| - |x|\big|(|fx|+|x|) + \big||fy| - |y|\big|(|fy|+|y|)\cr
& +\big||fx-fy| - |x-y|\big|(|fx-fy| + |x-y|)\\
&\le \varepsilon(2|x| + \varepsilon) + \varepsilon(2|y| + \varepsilon) + \varepsilon(2|x-y| +
\varepsilon)\\
&\le \varepsilon(4|x| + 4|y| + 3\varepsilon), 
\end{align*}
and the lemma follows.   $\square$

\alku
\label{iso}
\lem
Suppose that $0 \in Y \sub E	$ and that $\cspa Y = E$.  Suppose also that
$\varphi\colon Y
\to F$ is an iso\-metry with $\varphi(0)=0$. Then there is a unique extension of
$\varphi$ to an iso\-metry $T\colon E \to F$. Moreover, $T$ is linear.
\loppu

\pf
This result is well known. For the finite-dimensional case, see \cite[p. 9]{Re}. The
general case follows easily from this by considering restrictions of $\varphi$ to
finite subsets of $Y$ and making use of the uniqueness of $T$.  $\square$

\alku
\label{4.7}
\lau
Suppose that $A \sub E$ is an unbounded set such that $0\in A$ and $\mu(A) \ge 1/c'$.
Suppose also that $f\colon A \to F$ is an \eni with $f(0)=0$. Then there is a linear
iso\-metry $U\colon E \to F$ with $d(U,Pf) \le 2c'\varepsilon$, where $P\colon F \to
UE$ is the orthogonal projection.
\loppu

\proof
Set $X = \cd$ and let $u \in X$. Let $(x_j)$ be a sequence in $A$ converging
directionally to $u$. By \ref{4.5}, the sequence $(pfx_j)$ converges to a point $u' \in
S(1)$, and $u'$ is independent of the choice of $(x_j)$. Setting $\varphi u = u'$ we
thus obtain a well defined map $\varphi\colon X \to S(1)$. We show that $\varphi$ is
an iso\-metry.

Since $f(0)=0$ and since $f$ is an \eni, we have
$
\big| |fx|/|x| - 1\big| \le \varepsilon/|x|
$
for all $x \in A \sm\{  0\}$. Hence $|fx_j|/|x_j| \to 1$ whenever $x_j \in A$ and
$|x_j| \to\iy$. Consequently,
\beq
\label{4.a}
\varphi u = \lim_{j\to\iy} \frac{fx_j}{|x_j|}
\end{equation}
for each sequence $(x_j)$ in $A$ converging directionally to $u$.

Let $u,v \in X$ and choose sequences $(x_j), (y_j)$ in $A$ converging directionally to
$u$ and $v$, respectively. By
\ref{M} we have
\[
|fx_j \cdot fy_j - x_j \cdot y_j| \le 2\varepsilon(|x_j| + |y_j| + \varepsilon).
\]
Dividing by $|x_j||y_j|$ and letting $j\to\iy$ yields $\varphi u \cdot \varphi v = u
\cdot v$. Since the vectors $u,v,\varphi u,\varphi v$ have norm one, this implies that
$|\varphi u-\varphi v| = |u-v|$, and thus
$\varphi$ is an iso\-metry.

Extend $\varphi$ to an iso\-metry $\varphi_0\colon  X \cup \{  0\} \to $ by
$\varphi_0(0)=0$. Since
$\mu_1(X) = \mu(A) > 0$, we have $\cspa X = E$. By \ref{iso}, the map $\varphi_0$
extends to a linear iso\-metry $U \colon E\to F$.  It
remains to show that
$d(U,Pf) \le 2c'\varepsilon$.

Let $U_1\colon E \to UE$ be the bijective linear iso\-metry defined by $U$. 
Replacing $f$ by $f U_1^{-1}\colon UE \to F$ we may assume that $E \sub F$ and $U =
\id$. Let
$x
\in E$ and set
$\alpha = |Pfx-x|$. We must show that $\alpha \le 2c'\varepsilon$.

Choose a unit vector $e$ with $Pfx-x = \alpha e$ and let $\lambda > 1$. By the
definition \ref{2.2} of $\mu(A)$, there is $u \in X$ with $|u \cdot e| \ge
\mu(A)/\lambda \ge 1/\lambda c'$. Choose a sequence $(x_j)$ in $A$ converging
directionally to $u$.
  Then $fx_j/|x_j| \to U u = u$. By \ref{M} we have
\[
|fx \cdot fx_j - x \cdot x_j| \le 2\varepsilon (|x| + |x_j| + \varepsilon).
\]
Dividing by $|x_j| $ and letting $j\to\iy$ yields $|fx \cdot u - x \cdot u| \le
2\varepsilon$. Since $Pfx \cdot u = fx \cdot u$, this implies that
\[
\alpha/\lambda c'  \le \alpha |e \cdot u| = |(Pfx - x) \cdot u | \le
2\varepsilon.
\]
As $\lambda\to 1$, this gives $\alpha \le 2c'\varepsilon$.  $\square$

\alku
\label{4.7+1}
\cor
Let $f\colon A \to F$ be as in {\rm \ref{4.7}}. Then there is a linear map $S\colon F \to
E$ such that $|S|=1$ and $d(Sf,\id) \le 2c' \varepsilon$.
\loppu

\pf
Set $S = U_1^{-1}P$.  $\square$

\alku
\label{4.8}
Remark. \rm
In the case where $A=E$ and thus $\mu(A)=1$, Corollary \ref{4.7+1} was proved by
\cite[Th. 8]{Qi}
 with the bound $6\varepsilon$. A direct proof for this case with the bound
$2\varepsilon$ is given in Section 5; see \ref{5.3+1}. 
\loppu

\alku
\label{4.9}
\cor
Suppose that $A \sub \rn$ is an unbounded set such that $0\in A$ and $\mu(A) \ge
1/c'$. Suppose also that $f\colon A \to \rn$ is an \eni with $f(0)=0$. Then there is a
linear iso\-metry $U\colon \rn \to \rn$ with $d(U,f) \le 2c'\varepsilon$.  $\square$
\loppu

\alku
\label{4.10}
The Jung constant. \rm
The Jung constant $J(V)$ of a normed space $V$ is the infimum of all $r > 0$ such that
every set $Q \sub V$ with $d(Q) \le 2$ is contained in a ball of radius $r$. We have
always $1 \le J(V) \le 2$, and $J(\rn) = \sqrt{2n/(n+1)}$ by the classical result of
H.W.E. Jung \cite{Ju}. Furthermore, $J(E) = \sqrt 2$ for infinite-dimensional Hilbert
spaces $E$; see \cite[Th. 2]{Da} or \cite[p. 704]{Se}. 
\loppu

\alku
\label{4.11}
\lau
Suppose that $A \sub E$ is an unbounded set with $\mu(A) \ge 1/c' > 0$. Let $f\colon A
\to F$ be an \eni. Then there is an iso\-metry $T\colon E\to F$ onto a closed linear
subspace of $F$ such that
\[
d(T,Pf) \le J(F)c'\varepsilon \le \sqrt2 c'\varepsilon,
\]
where $P\colon F \to TE$ is the orthogonal projection.
\loppu

\proof
We follow the idea  of \cite[1.2]{Se}.  For each $a \in A$ we define $g_a\colon A-a
\to F$ by $g_a x = f(x+a) - fa$. Then $0 \in A-a$ and $g_a(0)=0$. Let $U_a\colon
E\to F$ be the linear iso\-metry given by \ref{4.7} for the map $g_a$.  

Let $u \in \cd$, and let $(x_j)$ be a sequence in $A$ converging directionally to $u$.
Then $u \in {\rm cd}\, (A-a)$ and $(x_j-a)$ converges directionally to $u$. By
(\ref{4.a}) we have
\[
U_au = \lim_{j\to\iy} \frac{g_a(x_j-a)}{|x_j-a|} = \lim_{j\to\iy}
\frac{fx_j-fa}{|x_j-a|} = \lim_{j\to\iy} \frac{fx_j}{|x_j|},
\]
and hence $U_au$ is independent of $a$ for each $u \in \cd$. Since $\cspa \cd = E$, it
follows that $U_a$ is independent of $a$, and we write $U=U_a$. 

Let $P\colon F \to UE$ be the orthogonal projection and set 
 $h = f-U\colon  A \to F$. Then $d(U,Pg_a) \le 2c'\varepsilon$ for each $a \in A$. For
all $a,b \in A$ we have
\[
|Pha-Phb| = |Pg_b(a-b) - U_b(a-b)| \le 2c'\varepsilon.
\]
Hence $d(PhA) \le 2c'\varepsilon$. Let $\lambda > 1$.  Then
$PhA$ is contained in the ball $\bar B(w,\lambda J(F)c'\varepsilon)$ for some $w \in UE$.
Setting
$Tx = Ux + w$ we obtain an isometry
$T\colon E\to F$ with $TE = UE$. For each $x \in A$ we have
\[
|Tx-Pfx| = |w-Phx| \le \lambda J(F)c'\varepsilon.
\]
As $\lambda\to 1$, this gives the theorem.  $\square$ 

\alku
\label{4.12}
\cor
Part $(2)$ \imp $(1)$ of Theorem {\rm \ref{main}} is true.  $\square$
\loppu

\section{A short proof for the Hyers-Ulam theorem}

The proof of the original Hyers-Ulam theorem in \cite{HU} is elementary but rather long.
We next give a considerably shorter elementary proof, which leads to the optimal
constant and also gives  the finite-dimensional version due to Bhatia-\v Semrl
\cite[Th. 1]{BS}. A crucial tool is our inequality \ref{M}. Otherwise, the proof is
self-contained.

In \ref{A.4} we give an improved version of the Hyers-Ulam theorem with relaxed
surjectivity condition.

\alku
\label{A.1}
\lau
Suppose that $E$ is a Banach space, that $F$ is a  Hilbert space, and that $f\colon E
\to F$ is an
\eni with
$f(0)=0$.  Then there is a
 linear isometry $T\colon E \to F$ such that $d(T,Pf) \le 2\varepsilon$, where $P\colon 
F \to TE$ is the orthogonal projection. Hence $E$ is isomorphic to a Hilbert space.

If $f$ is surjective or if $\dim E = \dim F < \iy$, then $T$ is surjective and $d(T,f)
\le 2\varepsilon$.
\loppu

\pf
We first show that the limit $Tx = \lim_{s\to\iy} f(sx)/s$ exists for each $x \in E$.
 Since $F$ is complete, it
suffices to show that
\beq
\label{delta}
\sup \{ | f(sx)/s - f(tx)/t|: t \ge s\} \to 0
\end{equation}
as $s \to\iy$. Let $0 < s \le t$. By \ref{M} we obtain
\begin{align*}
&|f(sx)/s-f(tx)/t|^2 = |f(sx)|^2/s^2 + |f(tx)|^2/t^2 - 2f(sx) \cdot f(tx)/st\\
& \le (s|x|+\varepsilon)^2/s^2 + (t|x|+\varepsilon)^2/t^2 -2|x|^2 +
4\varepsilon(|x|/t+|x|/s+\varepsilon/st)\\
& \le 12\varepsilon|x|/s + 6\varepsilon^2/s^2,
\end{align*}
which implies (\ref{delta}). Observe that $T(0)=0$.

We next show that the map $T\colon E \to F$ is an iso\-metry and hence linear. Let $x,y
\in E$. Since $f$ is an \eni,  we have 
\[
\big| |f(sx)-f(sy)| - |sx-sy|\big| \le \varepsilon.
\]
Dividing by $s$ and letting $s \to\iy$ yields $|Tx-Ty| = |x-y|$.

Let $T_1\colon E \to TE$ be the bijective linear iso\-metry defined by $T$.
 Replacing $f$ by $fT_1^{-1}\colon
TE \to F$ we may assume that $E \sub F$ and that $T = \id$. 
Let $P\colon F \to E$ be the orthogonal projection and let $x \in E$. We must show that 
\beq
\label{A.2}
|Pfx - x| \le 2\varepsilon.
\end{equation}

Set $\alpha = |Pfx-x|$ and choose a unit vector $u \in E$ with $\alpha u = Pfx-x$. By
\ref{M} we have
\[
|fx \cdot f(su) - x \cdot su| \le 2\varepsilon(|x| + s + \varepsilon)
\]
for all $s > 0$. Since $f(su)/s \to Tu = u$ as $s \to\iy$, this yields $|fx \cdot u -
x \cdot u| \le 2\varepsilon$. Since $fx \cdot u = Pfx \cdot u$, we obtain
\[
\alpha = \alpha u \cdot u = fx \cdot u - x \cdot u\le 2\varepsilon,
\]
and (\ref{A.2}) follows.

If $\dim E = \dim F < \iy$, then $T$ is surjective. Assume that $f$ is surjective and
that $E \ne F$. Choose a unit vector
$e
\in E^{\perp}$ and then $x \in E $ with $fx = (3\varepsilon + 1)e$. 
Since $f$ is an \eni and since $Pfx = 0$, we obtain by (\ref{A.2}) the contradiction
\[
3\varepsilon + 1 = |fx| \le |x| + \varepsilon \le 3\varepsilon. \quad \square
\]

\alku
\label{5.3+1}
 Remark. \rm
Setting $S = T_1^{-1} P$ in \ref{A.1} we obtain a linear map $S\colon F \to E$ with $|S|
= 1$ and $d(Sf,\id) \le 2\varepsilon$.
\loppu

We finally show that the surjectivity condition of \ref{A.1} can be replaced by the
weaker condition $\tau(fE) < 1$, where $\tau$ was defined in \ref{alt}. Observe that
$\tau(F) = 0$.

\alku
\label{A.4}
\lau
Suppose that $E$ is a Banach space, that   $F$ is a  Hilbert space, and that $f\colon
E
\to F$ is an
\eni such that $f(0)=0$ and $\tau(fE) < 1$. Then there is a surjective linear isometry
$T\colon E
\to F$ with $d(T,f) \le 2\varepsilon$.
\loppu

\pf
The proof of \ref{A.1} is valid until the last paragraph. Assume again that $E \ne F$,
and choose a unit vector $e \in E^\perp$. Choose a number $q$ with $\tau(fE) < q < 1$.
Since $\liminf_{|t| \to\iy} d(te,fE)/|t| < q$, there are sequences $(t_j)$ in ${\sf R}$
and $(x_j)$ in $E$ such that $|t_j| \to\iy$ and 
\beq
\label{A.5}
|t_je-fx_j| \le q|t_j|
\end{equation}
for all $j$. Setting $y_j = fx_j$ we have
\beq
\label{A.6}
(1-q)|t_j| \le |y_j| \le (1+q)|t_j| \le 2|t_j|
\end{equation}
for all $j$. 

Since $|Py_j - x_j| \le 2\varepsilon$ by (\ref{A.2}), we have $|Py_j| \ge |x_j| -
2\varepsilon \ge |y_j| - 3\varepsilon$. Since $|y_j| \to\iy$ by (\ref{A.6}), we have
$|y_j| \ge 3\varepsilon$
 for large $j$, and then $|Py_j|^2 \ge |y_j|^2 - 6\varepsilon|y_j|$. Since $|y_j|^2 =
|Py_j|^2 + |y_j - Py_j|^2$, we get
\[
e \cdot y_j = e \cdot (y_j-Py_j) \le |y_j-Py_j| \le \sqrt{6\varepsilon|y_j|} \le 4
\sqrt{\varepsilon t_j}
\]
by (\ref{A.6}). Since (\ref{A.5}) implies that
\[
t_j^2 + |y_j|^2 -2t_j e \cdot y_j \le q^2 t_j^2,
\]
we obtain $(1-q^2)t_j^2 \le 8t_j\sqrt{\varepsilon t_j}$. Dividing by $t_j^2$ and letting
$j \to\iy$ gives the desired contradiction.  $\square$

\alku
\label{A.7}
\cor
If $f\colon E \to F$ is a near\-isometry, then $\tau(fE) \in \{  0,1\}$.
\loppu

\pf
Since $\tau$ is invariant under translations, we may assume that $f(0)=0$. If
$\tau(fE) < 1$, it follows from \ref{A.4} that $d(y,fE) \le 2\varepsilon$ for all $y \in
F$, and hence $\tau(fE) = 0$.  $\square$

\section{A remark on Banach spaces}

The proofs in the preceding sections make substantial use of the inner product, and
there seems to be no easy way to extend them for general Banach spaces. However, we show
that instead of considering surjective near\-isometries between Banach spaces $E,F$, it
suffices to consider maps $f\colon A \to F$ such that the sets $E \sm A$ and $F \sm fA$
are bounded. This means that the theorem of Hyers-Ulam-Gevirtz-Omladi\v c-\v Semrl
\cite[15.2]{BL} for Banach spaces is not, after all, a global result but a local
property of maps near the point $\iy$. The result was suggested to the author by O.
Martio.

\alku
\label{6.1}
\lau
Suppose that $E$ and $F$ are Banach spaces, that $A \sub E$, 
 and that $f\colon A \to F$ is an $\varepsilon$-near\-isometry such
that the sets $E \sm A$ and $F \sm fA$ are bounded. Then there is a surjective iso\-metry
$T\colon E \to F$ with $d(T,f) \le 2\varepsilon$. For each $x_0\in A$ we can choose $T$ so that
$Tx_0 = fx_0$.
\loppu

\pf
We may assume that $x_0 = 0,\ fx_0 = 0$. Choose a number $R > \varepsilon$ such that $E \sm A
\sub B(R)$ and $F \sm fA \sub B(2R)$. Define $f_1\colon E \to F$ by $f_1x = 2x$ for $|x| \le R$
and by
$f_1x = fx$ for $|x| > R$. Then $f_1$ is a surjective $K$-near\-isometry with $K =
3R+\varepsilon$. By \cite[15.2]{BL}, there is a linear surjective isometry
$T\colon E \to F$ with $d(T,f_1) \le 2K$. Hence $d(T,f) < \iy$. Now an easy
modification of the proof of
\cite[15.2]{BL} (see
\cite[p. 362]{BL}) shows that $d(T,f) \le 2\varepsilon$.  $\square$ 

\medskip
If $E$ and $F$ are Banach spaces of the same finite dimension, then every \eni $f\colon
E \to F$ with $f(0)=0$ can be approximated by a linear isometry $T\colon E \to F$ with
$d(T,f) \le 2\varepsilon$ \cite[Th. 1]{Di}.

\alku
\label{6.2}
Open problem. \rm
Suppose that $E$ and $F$ are finite-dimensional Banach spaces with $\dim E = \dim F <
\iy$, that
$A
\sub E$ is a half space and that $f\colon A \to F$ is an \eni. Does there exist an
isometry
$T\colon E \to F$ such that $d(T,f) \le K\varepsilon$ with some universal constant $K$?
\loppu

\small
\renewcommand{\baselinestretch}{0.9}

\noindent Preprints of the author can be downloaded
from www.helsinki.fi/$^\sim$jvaisala/preprints.html.
\bigskip

\noindent Matematiikan laitos\\ 
Helsingin yliopisto\\ 
PL 4, Yliopistonkatu 5\\ 
00014 Helsinki, Finland\\
\texttt{jvaisala@cc.helsinki.fi}


\begin{thebibliography}{xxx} 



\bibitem [ATV]{ATV} P. Alestalo, D.A. Trotsenko and J. V\"ais\"al\"a,
\textit{Isometric approximation}, Israel J. Math. \textbf{125} (2001), 61--82.

\bibitem [BL]{BL} Y. Benyamini and J. Lindenstrauss, \textit{Geometric nonlinear
functional analysis I}, AMS Colloquium Publications 48, 2000.

\bibitem [B\v S]{BS}  R. Bhatia and P. \v Semrl, \textit{Approximate isometries on
Euclidean spaces}, Amer. Math. Monthly \textbf{104} (1997), 497--504.

\bibitem [Da]{Da} J. Dane\v s,  \textit{On the radius of a set in a Hilbert space},
Comment. Math. Univ. Carolin. \textbf{25} (1984), 355--362.

\bibitem [Di]{Di} S.J. Dilworth, \textit{Approximate isometries in finite-dimensional
normed spaces}, Bull. London Math. Soc. \textbf{31} (1999), 471--476.

\bibitem [HV]{HV} T. Huuskonen and J. V\"ais\"al\"a, \textit{Hyers-Ulam constants of
Hilbert spaces}, preprint.
      
\bibitem [HU]{HU} D.H. Hyers and S.M. Ulam, \textit{On approximate isometries}, Bull.
Amer. Math. Soc. \textbf{51} (1945), 288--292.

\bibitem [Jo]{Jo} F. John, \textit{Rotation and strain}, Comm. Pure Appl. Math.
\textbf{14} (1961), 391--413.

\bibitem [Ju]{Ju} H.W.E. Jung, \textit{\"Uber die kleinste Kugel, die eine r\"aumliche
Figur einschliesst}, J. Reine Angew. Math. \textbf{123} (1901), 241--257.

\bibitem [Ma]{Ma} E. Matou\v skov\'a, \textit{Almost isometries of balls}, preprint.

\bibitem [Qi]{Qi} S. Qian, \textit{$\varepsilon$-isometric embeddings}, Proc. Amer.
Math. Soc. \textbf{123} (1995), 1797--1803.

\bibitem [Re]{Re} E.G. Rees, \textit{Notes on geometry}, Springer, 1983.

\bibitem [\v Se]{Se} P. \v Semrl, \textit{Hyers-Ulam stability of isometries}, Houston
J. Math. \textbf{24} (1998), 699--706. 

\bibitem [V\"a1]{solar} J. V\"ais\"al\"a, \textit{Isometric approximation property in
euclidean spaces},  Israel J. Math., to appear.

\bibitem [V\"a2]{survey} J. V\"ais\"al\"a, \textit{A survey of near\-isometries}, 
Report. Univ. Jyv\"askyl\"a, to appear.


\end{thebibliography}
\end{document}